\documentclass{amsart}
\usepackage{amsmath,amsfonts,latexsym,amssymb}

\newtheorem{theorem}{Theorem}
\newtheorem{lemma}{Lemma}

\theoremstyle{remark}

\newtheorem*{ack}{Acknowledgements}

\begin{document}

\markboth{Ritabrata Munshi}{Shifted convolution of divisor function $d_3$ and Ramanujan $\tau$ function}
\title[Shifted convolution of divisor function $d_3$ and Ramanujan $\tau$ function]{Shifted convolution of divisor function $d_3$\\ and Ramanujan $\tau$ function}

\author{Ritabrata Munshi}   
\address{School of Mathematics, Tata Institute of Fundamental Research, 1 Dr. Homi Bhabha Road, Colaba, Mumbai 400005, India.}     
\email{rmunshi@math.tifr.res.in}


\subjclass{11F66, 11M41}
\keywords{divisor functions, Ramanujan $\tau$ function, shifted convolution sum, circle method}

\maketitle


\section{Introduction}
\label{intro}

This note can be viewed as a bridge between the work of Pitt \cite{P} and my recent paper \cite{Mu0}. In \cite{P} Pitt considers the sum
$$
\Psi(f,x)=\sum_{1\leq n\leq x}d_3(n)a(rn-1) 
$$
where $d_3(n)$ is the divisor function of order $3$ (the coefficients of the Dirichlet series $\zeta(s)^3$), $a(m)$ is the normalized Fourier coefficients of a holomorphic cusp form $f$, and $r$ is a positive integer. Without loss one may take $a(m)=\tau(m)/m^{11/2}$ where $\tau$ is the Ramanujan function, as one does not expect any new complication to arise while dealing with Fourier coefficients of a general holomorphic cusp form. The trivial bound is given by $O(x^{1+\varepsilon})$. Pitt \cite{P} proved 
$$
\Psi(f,x)\ll x^{71/72+\varepsilon}
$$
where the implied constant is uniform with respect to $r$ in the range $0<r\ll X^{1/24}$. This sum is intrinsically related to the generalized Titchmarsh divisor problem, where one seeks to estimate the sum (see \cite{P}, \cite{P2})
$$
\sum_{\substack{p<x\\p\;\text{prime}}}a(p-1).
$$\\

In this paper we will use our method from \cite{Mu0}, \cite{Mu} to prove the following (improved bound).
\begin{theorem}
\label{mthm}
For $r\ll X^{1/10}$ we have
$$
\Psi(\Delta,x)\ll r^{\frac{2}{7}}X^{\frac{34}{35}+\varepsilon},
$$
where the implied constant depends only on $\varepsilon$.\\
\end{theorem}

\ack
I thank the organizers of `The Legacy of Srinivasa Ramanujan', in particular M.S. Raghunathan and Dipendra Prasad for their kind invitation.  \\


\section{Twisted Voronoi summation formulae}
\label{prelim}
 
Let 
$$
\Delta(z)=e(z)\prod_{n=1}^\infty \left(1-e(nz)\right)^{24}
$$
be the Ramanujan's $\Delta$ function. Here we are using the standard notation $e(z)=e^{2\pi iz}$. The function $\Delta(z)$ is a cusp form for $SL(2,\mathbb Z)$ of weight $12$. The Ramanujan $\tau$ function is defined as the Fourier coefficients of $\Delta(z)$, namely 
$$
\Delta(z)=\sum_{n=1}^\infty \tau(n)e(nz).
$$ 
Ramanujan conjectured, and later Deligne proved, that $|\tau(p)|\leq 2p^{11/2}$ for any prime number $p$. In light of this bound it is natural to define the normalized $\tau$ function as
$$
\tau_0(n)=\tau(n)/n^{11/2}.
$$
Using the modularity of $\Delta(z)$ one can establish the following Voronoi type summation formula for $\tau_0(n)$.  
\begin{lemma}
\label{voronoi2}
Let $q$ be a positive integer, and $a$ be an integer such that $(a,q)=1$. Let $g$ be a compactly supported smooth function on $\mathbb R_+$. We have
\begin{align}
\label{voronoi2-eqn}
\sum_{m=1}^\infty \tau_0(m)e_q\left(am\right)g(m)=\frac{2\pi}{q}\sum_{m=1}^\infty \tau_0(m)e_q\left(-\bar{a}m\right)G\left(\frac{m}{q^2}\right)
\end{align}
where $\bar{a}$ is the multiplicative inverse of $a\bmod{q}$, $e_q(z)=e(z/q)$ and
\begin{align*}
G(y)=\int_0^\infty g(x)J_{11}\left(4\pi\sqrt{xy}\right)dx.
\end{align*}
Here $J_{11}(z)$ is the Bessel functions in standard notations.
\end{lemma}

If $g$ is supported in $[AY,BY]$ (with $0<A<B$), satisfying $y^jg^{(j)}(y)\ll_j 1$, then the sum on the right hand side of \eqref{voronoi2} is essentially supported on $m\ll q^2(qY)^{\varepsilon}/Y$. The contribution from the tail $m\gg q^2(qY)^{\varepsilon}/Y$ is negligibly small. For smaller values of $m$ we will use the trivial bound $G(m/q^2)\ll Y$.\\

A similar Voronoi type summation formula for the divisor function $d_3(n)$ is also known (see Ivic \cite{I}). Let $f$ be a compactly supported smooth function on $\mathbb R_+$, and let $\tilde f(s)=\int_0^\infty f(x)x^{s}dx$. We define
\begin{align}
\label{gl}
F_{\pm}(y)=\frac{1}{2\pi i}\int_{\left(\frac{1}{8}\right)}(\pi^3 y)^{-s}\frac{\Gamma^3\left(\frac{1\pm 1+2s}{4}\right)}{\Gamma^3\left(\frac{3\pm 1-2s}{4}\right)}\tilde f(-s)ds.
\end{align}
\begin{lemma}
\label{voronoi3}
Let $f$ be a compactly supported smooth function on $\mathbb R_+$, we have
\begin{align}
\label{voronoi3-eqn}
\sum_{n=1}^\infty d_3(n)e_q\left(an\right)f(n)=&\frac{1}{q}\int_0^\infty P(\log y,q)f(y)\mathrm{d}y\\ 
\nonumber &+\frac{\pi^{3/2}}{2q^3}\sum_\pm\sum_{n=1}^\infty D_{3,\pm}(a,q;n)F_\pm\left(\frac{n}{q^3}\right),
\end{align}
where $P(y,q)=A_0(q)y^2+A_1(q)y+A_2(q)$ is a quadratic polynomial whose coefficients depend only on $q$ and satisfy the bound $|A_i(q)|\ll q^\varepsilon$. Also $D_{3,\pm}(a,q;n)$ are given by
$$
\sum_{n_1n_2n_3=n}\mathop{\sum\sum\sum}_{b,c,d=1}^q\left\{e_q(bn_1+cn_2+dn_3+abcd)\mp e_q(bn_1+cn_2+dn_3-abcd)\right\}.
$$
\end{lemma}

Suppose $f$ is supported in $[AX,BX]$, and $x^jf^{(j)}(x)\ll_j H^j$. Then the sums on the right hand side of \eqref{voronoi3-eqn} are essentially supported on $n\ll q^3H(qX)^{\varepsilon}/X$. The contribution from the tail $n\gg q^3H(qX)^{\varepsilon}/X$ is negligibly small. This follows by estimating the integral $F_{\pm}(y)$ by shifting the contour to the right. For smaller values of $n$ we shift the contour to left upto $\sigma=\varepsilon$.\\


\section{Setting up the circle method}
\label{cm}

As in \cite{Mu0}, we will be using a variant Jutila's version of the circle method. For any set $S \subset \mathbb R$, let $\mathbb I_S:\mathbb R\rightarrow \{0,1\}$ be defined by $\mathbb I_S(x)=1$ for $x\in S$ and $0$ otherwise. For any collection of positive integers $\mathcal Q \subset [1,Q]$ (which we call the set of moduli), and a positive real number $\delta$ in the range $Q^{-2}\ll \delta \ll Q^{-1} $, we define the function
$$
\tilde I_{\mathcal Q,\delta} (x)=\frac{1}{2\delta L}\sum_{q\in\mathcal Q}\;\sideset{}{^\star}\sum_{a\bmod{q}}\mathbb{I}_{[\frac{a}{q}-\delta,\frac{a}{q}+\delta]}(x),
$$
where $L=\sum_{q\in\mathcal Q}\phi(q)$. This is an approximation for $\mathbb I_{[0,1]}$ in the following sense (see \cite{J-1}): 

\begin{lemma}
\label{jutila-lemma}
We have
\begin{align}
\label{jutila}
\int_0^1\left|1-\tilde I_{\mathcal Q,\delta}(x)\right|^2dx\ll \frac{Q^{2+\varepsilon}}{\delta L^2}.
\end{align}
\end{lemma}

\vspace{.5cm}

Instead of studying the sum in Theorem \ref{mthm} we examine the related smoothed sum over dyadic segment
$$
D=\sum_{n=1}^\infty d_3(n)\tau_0(rn-1)W(n/X)
$$
where $W$ is a non-negative smooth function supported in $[1-H^{-1},2+H^{-1}]$ (we will chose $H=X^\theta$ optimally later), with $W(x)=1$ for $x\in [1,2]$ and satisfying $W^{(j)}(x)\ll_j H^j$. Clearly we have
\begin{align}
\label{prop}
\Psi(\Delta,x)=D+O(X^{1+\varepsilon}/H).
\end{align} 
In the rest of the paper we will prove a compatible bound for $D$.\\

Let $V$ be a smooth function supported in $[1/2,3]$ satisfying $V(x)=1$ for $x\in [3/4,5/2]$, $V^{(j)}(x)\ll_j 1$, and let $Y=rX$. Then we have
\begin{align*}
D=&\mathop{\sum\sum}_{n,m=1}^\infty d_3(n)\tau_0(m)W\left(\frac{n}{Y}\right)V\left(\frac{m}{Y}\right)\delta(rn-1,m)\\
=&\int_0^1e(-x)\left[\mathop{\sum}_{n=1}^\infty d_3(n)e(xrn)W\left(\frac{n}{X}\right)\right]\left[\mathop{\sum}_{m=1}^\infty \tau_0(m)e(-xm)V\left(\frac{m}{Y}\right)\right]dx.
\end{align*}
Let $\mathcal Q$, which we choose carefully later, be a collection of moduli of size $Q$. Suppose $|\mathcal Q|\gg Q^{1-\varepsilon}$, so that $L=\sum_{q\in\mathcal Q}\phi(q)\gg Q^{2-\varepsilon}$. Let $\delta=Y^{-1}$, and define
\begin{align}
\label{approx-0}
\tilde D:=\int_0^1\tilde I_{\mathcal Q,\delta}(x)e(-x)\left[\mathop{\sum}_{n=1}^\infty d_3(n)e(xrn)W\left(\frac{n}{X}\right)\right]\left[\mathop{\sum}_{m=1}^\infty \tau_0(m)e(-xm)V\left(\frac{m}{Y}\right)\right]dx.
\end{align}
It follows that
$$
\tilde D=\frac{1}{2\delta}\int_{-\delta}^\delta \tilde D(\alpha)e(-\alpha)\mathrm{d}\alpha,
$$
where
\begin{align}
\label{approx}
\tilde D(\alpha)=\frac{1}{L}\sum_{q\in\mathcal Q}\;\sideset{}{^\star}\sum_{a\bmod{q}}e_q(-a)&\left[\mathop{\sum}_{n=1}^\infty d_3(n)e_q(arn)e(\alpha rn)W\left(\frac{n}{X}\right)\right]\\
\nonumber &\times \left[\mathop{\sum}_{m=1}^\infty \tau_0(m)e_q(-am)e(-\alpha m)V\left(\frac{m}{Y}\right)\right].
\end{align}\\

In circle method we approximate $D$ by $\tilde D$, and then try to estimate the latter sum. Lemma~\ref{jutila-lemma} gives a way to estimate the error in this process. More precisely we have
\begin{align}
\label{error}
\left|D-\tilde D\right|\ll \int_0^1\left|\mathop{\sum}_{n=1}^\infty d_3(n)e(xrn)W\left(\frac{n}{X}\right)\right|\left|\mathop{\sum}_{m=1}^\infty \tau_0(m)e(-xm)V\left(\frac{m}{Y}\right)\right|\left|1-\tilde I_{\mathcal Q, \delta}(x)\right|dx,
\end{align}
Using the well-known point-wise uniform bound 
$$
\mathop{\sum}_{m=1}^\infty\tau_0(m)e(-xm)V\left(\frac{m}{Y}\right)\ll Y^{\frac{1}{2}+\varepsilon}
$$
it follows that the right hand side of \eqref{error} is bounded by
\begin{align*}
\ll Y^{\frac{1}{2}+\varepsilon}\int_0^1\left|\mathop{\sum}_{n=1}^\infty d_3(n)e(xrn)W\left(\frac{n}{X}\right)\right|\left|1-\tilde I(x)\right|dx.
\end{align*}
Now we apply Cauchy and Lemma \ref{jutila-lemma}, with $Q=YX^{-\frac{1}{2}+\delta}$ for any $\delta>0$, to arrive at the following:
\begin{lemma}
\label{first-lemma}
We have
\begin{align}
\label{exp1}
D=\tilde D+O\left(X^{1-\delta+\varepsilon}\right).
\end{align}
\end{lemma}
We can now decide what will be the optimal choice for $H$. Naturally we wish to take $H$ as small as possible to aid in our analysis of $\tilde D$. Matching the error term in Lemma~\ref{first-lemma} and that in \eqref{prop} we pick $H=X^{\delta}$.


\section{Estimation of $\tilde D$}
\label{gl1_twists}

Now we apply Voronoi summations on the sums over $m$ and $n$. This process gives rise to several terms as noted in Section \ref{prelim} - Lemma \ref{voronoi2} and Lemma \ref{voronoi3}. As far as our analysis is concerned we can focus our attention on two such term, namely  
\begin{align}
\label{approx-pw-0}
\tilde D_0(\alpha)=\frac{2\pi}{L}\sum_{q\in\mathcal Q}\frac{1}{q^2}\sum_{m=1}^\infty \tau_0(m)S(1,m;q)G\left(\frac{m}{q^2}\right)\int_0^\infty P(\log x)f(x)\mathrm{d}x,
\end{align}
which is the zero frequency contribution ($S(1,m;q)$ is the Kloosterman sum), and 
\begin{align}
\label{approx-pw-1}
\tilde D_1(\alpha)=\frac{\pi^{5/2}}{L}\sum_{q\in\mathcal Q}&\frac{1}{q^4}\sum_{m=1}^\infty \tau_0(m)\sum_{n=1}^\infty \mathcal S^\star(m,n;q)G\left(\frac{m}{q^2}\right)F_+\left(\frac{n}{q^3}\right),
\end{align}
where the character sum is given by
\begin{align}
\label{main-char}
\mathcal S^\star(m,n;q):=\sideset{}{^\star}\sum_{a\bmod{q}}e_q(-a-\bar{a}m)
\sum_{n_1n_2n_3=n}\mathop{\sum\sum\sum}_{b,c,d=1}^q e_q(bn_1+cn_2+dn_3+abcdr)
\end{align}
Also here we are taking
$$
g(y)=V\left(\frac{y}{Y}\right)e(-\alpha y),\;\;\;\text{and}\;\;\;f(x)=W\left(\frac{x}{X}\right)e(-\alpha rx).
$$
The functions $G$ and $F_+$ are defined in Lemma \ref{voronoi2} and Lemma \ref{voronoi3} respectively. It follows that in both the sums \eqref{approx-pw-0} and \eqref{approx-pw-1}, the sum over $m$ essentially ranges upto $m\ll Q^2Y^{-1+\varepsilon}=YX^{-1+2\delta+\varepsilon}$. The tail contribution is negligibly small. So using the Weil bound for the Kloosterman sums it follows that
$$
\tilde D_0(\alpha)\ll X^{1+\varepsilon}/\sqrt{Q} 
$$
which is smaller than the bound in Lemma~\ref{first-lemma} (as $Q=YX^{-\frac{1}{2}+\delta}>X^{2\delta}$ or $Y>X^{\frac{1}{2}+\delta}$). One can use Deligne's theory to show that there is square root cancellation in the character sum \eqref{main-char}. But this is not enough to establish a satisfactory bound for $\tilde{D}_1(\alpha)$. \\

Following \cite{Mu0} we will now make an appropriate choice for the set of moduli. We choose  $\mathcal Q$ to be the product set $\mathcal Q_1 \mathcal Q_2$, where $\mathcal Q_i$ consists of primes in the dyadic segment $[Q_i,2Q_i]$ (and not dividing $r$) for $i=1,2$, and $Q_1Q_2=Q$. Also we pick $Q_1$ and $Q_2$ (whose optimal sizes will be determined later) so that the collections $\mathcal Q_1$ and $\mathcal Q_2$ are disjoint. \\

Suppose $q=q_1q_2$ with $q_i\in\mathcal Q_i$. The character sum $\mathcal S^\star(m,n;q)$ splits as a product of two character sums with prime moduli. The one modulo $q_1$ looks like (after a change of variables)
$$
\mathcal S^\dagger(m,n,q_2;q_1)=\sideset{}{^\star}\sum_{a\bmod{q_1}}e_{q_1}(-\bar{q_2}^3a-q_2\bar{a}m)
\sum_{n_1n_2n_3=n}\mathop{\sum\sum\sum}_{b,c,d=1}^{q_1} e_{q_1}(bn_1+cn_2+dn_3+abcdr).
$$
Now let us consider the case where $q_1|n$. Suppose $q_1|n_1$, then summing over $b$ we arrive at
$$
q_1\mathop{\sum}_{d=1}^{q_1} e_{q_1}(dn_3)+q_1\mathop{\sum}_{c=1}^{q_1} e_{q_1}(cn_2)-q_1.
$$ 
This sum is bounded by $q_1(q_1,n_2n_3)$. Then using Weil bound for Kloosterman sums we conclude that
$$
\mathcal S^\dagger(m,q_1n,q_2;q_1)\ll q_1^{3/2}(q_1,n)d_3(n).
$$
On the other hand if $q_1\nmid n$ then we arrive at the following expression for the character sum after a change of variables 
$$
\mathcal S^\dagger(m,n,q_2;q_1)=d_3(n)\sideset{}{^\star}\sum_{a\bmod{q_1}}e_{q_1}(-\bar{q_2}^3a-q_2\bar{a}m)\mathop{\sum\sum\sum}_{b,c,d=1}^{q_1} e_{q_1}(b+c+d+\bar{n}abcdr).
$$
Summing over $b$ we arrive at
$$
\mathcal S^\dagger(m,n,q_2;q_1)=d_3(n)q_1\sideset{}{^\star}\sum_{a\bmod{q_1}}e_{q_1}(-\bar{q_2}^3a-q_2\bar{a}m)S(1,-n\overline{ar};q_1).
$$
This can be compared with the character sums which appear in \cite{Mu0} and \cite{P}. Strong bounds (square root cancellation) have been established for this sums using Deligne's result. In the light of this, it follows that to estimate the contribution of those $n$ in \eqref{approx-pw-1} with $(n,q)\neq 1$ it is enough to look at the sum
$$
\frac{1}{L}\sum_{q\in\mathcal Q}\frac{1}{Q^4}\sum_{0<m\ll Q^2/Y}\;\sum_{0<n\ll Q^3H/\min\{Q_1,Q_2\}X} Q^{3/2}\sqrt{\max\{Q_1,Q_2\}}XY\sqrt{H}.
$$ 
The last sum is bounded by $O(Q^{2+\varepsilon} H^{3/2}\min\{Q_1,Q_2\}^{-3/2})$, and we get
\begin{align*}
\tilde D_1(\alpha)=\frac{\pi^{5/2}}{L}\sum_{q\in\mathcal Q}\frac{1}{q^3}\sum_{m=1}^M &\tau_0(m)\sum_{\substack{n=1\\(n,q)=1}}^N d_3(n)\mathcal S^\dagger(m,n;q)G\left(\frac{m}{q^2}\right)F_+\left(\frac{n}{q^3}\right)+O\left(\frac{r^2X^{1+7\delta/2+\varepsilon}}{\min\{Q_1,Q_2\}^{3/2}}\right),
\end{align*}
where 
$$
\mathcal{S}^\dagger (m,n;q)=\sideset{}{^\star}\sum_{a\bmod{q}}e_{q}(-a-\bar{a}m)S(1,-n\overline{ar};q),
$$
$M=Q^{2+\varepsilon}Y^{-1}=rX^{2\delta+\varepsilon}$ and $N=Q^{3+\varepsilon}HX^{-1}=r^3X^{1/2+4\delta+\varepsilon}$. Next we observe that we can now remove the coprimality restriction $(n,q)=1$, without worsening the error term. Here we are using square root cancellation in the character sum $\mathcal{S}^\dagger(m,n;q)$. We get
\begin{align}
\label{approx-pw-11}
\tilde D_1(\alpha)=\frac{\pi^{5/2}}{L}\sum_{q\in\mathcal Q}\frac{1}{q^3}\sum_{m=1}^M &\tau_0(m)\sum_{n=1}^N d_3(n)\mathcal S^\dagger(m,n;q)G\left(\frac{m}{q^2}\right)F_+\left(\frac{n}{q^3}\right)\\
\nonumber &+O\left(\frac{r^2X^{1+7\delta/2+\varepsilon}}{\min\{Q_1,Q_2\}^{3/2}}\right),
\end{align}\\


\section{Estimation of $\tilde D_1(\alpha)$ : Final analysis}

Applying Deligne's bound for $\tau(m)$, the problem now reduces to estimating 
\begin{align*}
\frac{1}{Q^5}\sum_{q_2\in\mathcal Q_2}\sum_{1\leq m\leq M} \mathop{\sum}_{1\leq n\leq N} \left|\sum_{q_1\in\mathcal Q_1}\mathcal S^\dagger(m,n;q)G\left(\frac{m}{q^2}\right)F_+\left(\frac{n}{q^3}\right)\right|,
\end{align*}
where $q=q_1q_2$. Applying Cauchy inequality we get
\begin{align}
\label{dh}
\tilde D_1(\alpha)\ll \frac{\sqrt{N}}{Q^5}\sum_{q_2\in\mathcal Q_2}\sum_{1\leq m\leq M}\tilde D^\sharp(m,q_2)^{\frac{1}{2}},
\end{align}
where
\begin{align*}
\tilde D^\sharp(m,q_2)=\mathop{\sum}_{n\in \mathbb Z}h\left(n\right)\left|\sum_{q_1\in\mathcal Q_1}\mathcal S^\dagger(m,n;q_1q_2)G\left(\frac{m}{q_1^2q_2^2}\right)F_{+}\left(\frac{n}{q_1^3q_2^3}\right)\right|^2.
\end{align*}
Here $h$ a is non-negative smooth function on $(0,\infty)$, supported on $[1/2,2N]$, and such that $h(x)=1$ for $x\in [1,N]$ and $x^jh^{(j)}(x)\ll 1$.\\ 

Opening the absolute square and interchanging the order of summations we get
\begin{align*}
\tilde D^\sharp(m,q_2)=&\sum_{q_1\in\mathcal Q_1}\sum_{\tilde q_1\in\mathcal Q_1}G\left(\frac{m}{q_1^2q_2^2}\right)\bar G\left(\frac{m}{\tilde q_1^2q_2^2}\right)\\
&\times \mathop{\sum}_{n\in \mathbb Z}h\left(n\right)\mathcal S^\dagger(m,n;q_1q_2)\bar{\mathcal S^\dagger}(m,n;\tilde q_1q_2) F_+\left(\frac{n}{q_1^3q_2^3}\right)\bar F_+\left(\frac{n}{\tilde q_1^3q_2^3}\right).
\end{align*}
Applying Poisson summation on the sum over $n$ with modulus $q_1\tilde q_1q_2$, we get 
\begin{align}
\label{poisson-last}
\frac{1}{q_2}\sum_{q_1\in\mathcal Q_1}&\sum_{\tilde q_1\in\mathcal Q_1}\frac{1}{q_1\tilde q_1}G\left(\frac{m}{q_1^2q_2^2}\right)\bar G\left(\frac{m}{\tilde q_1^2q_2^2}\right)\sum_{n\in\mathbb Z}\mathcal T(m,n;q_1,\tilde q_1,q_2)\mathcal I(n;q_1,\tilde q_1,q_2).
\end{align}
The character sum is given by
$$
\mathcal T(m,n;q_1,\tilde q_1,q_2)=\sum_{\alpha\bmod{q_1\tilde q_1q_2}}\mathcal S^\dagger(m,\alpha;q_1q_2)\bar{\mathcal S^\dagger}(m,\alpha;\tilde q_1q_2)e_{q_1\tilde q_1q_2}(n\alpha),
$$
and the integral is given by
$$
\mathcal I(n;q_1,\tilde q_1,q_2)=\int_{\mathbb R} h\left(x\right)F_+\left(\frac{x}{q_1^3q_2^3}\right)\bar F_+\left(\frac{x}{\tilde q_1^3q_2^3}\right)e_{q_1\tilde q_1q_2}(-nx)dx.
$$
Integrating by parts repeatedly one shows that the integral is negligibly small for large values of $|n|$, say $|n|\geq X^{2013}$. Observe that differentiating under the integral sign in \eqref{gl}, one can show that $y^jF_+^{(j)}(y)\ll_j XH$. So we have the bound
$$
\mathcal I(n;q_1,\tilde q_1,q_2)\ll\frac{X^2HQ^3}{|n|}.
$$
The following lemma now follows from \eqref{poisson-last}.
\begin{lemma}
\label{rem-sum}
We have
\begin{align*}
\tilde D^\sharp(m,q_2)\ll (XY)^2H\sum_{q_1\in\mathcal Q_1}&\sum_{\tilde q_1\in\mathcal Q_1}\left\{\sum_{1\leq |n|\leq X^{2013}}\frac{H}{|n|}|\mathcal T(m,n;q_1,\tilde q_1,q_2)|+\frac{N}{QQ_1}|\mathcal{T}(m,0;q_1,\tilde{q}_1,q_2)|\right\}\\
&+X^{-2013}.
\end{align*}
\end{lemma}

It now remains to estimate the character sum. This has been done in \cite{Mu0}. We summarize the result in the following lemma.\\

\begin{lemma}
\label{char-sum-1}
For $q_1\neq \tilde q_1$, the character sum $\mathcal T(m,n;q_1,\tilde q_1,q_2)$ vanishes unless $(n,q_1\tilde q_1)=1$, in which case we have
$$
\mathcal T(m,n;q_1,\tilde q_1,q_2)\ll q_1^{\frac{3}{2}}\tilde q_1^{\frac{3}{2}}q_2^{\frac{5}{2}}(n,q_2)^{\frac{1}{2}}.
$$
The character sum $\mathcal T(m,n;q_1,q_1,q_2)$ vanishes unless $q_1|n$, in which case we have
$$
\mathcal T(m,q_1n';q_1,q_1,q_2)\ll q_1^{\frac{5}{2}}q_2^{\frac{5}{2}}\sqrt{(n',q_1q_2)}.
$$\\
\end{lemma}

It follows from Lemma \ref{char-sum-1}, that
\begin{align}
\label{bd-char-1}
\mathop{\sum_{q_1\in\mathcal Q_1}\sum_{\tilde q_1\in\mathcal Q_1}}_{q_1\neq \tilde q_1}\sum_{1\leq |n|\leq X^{2013}}|\frac{\mathcal T(m,n;q_1,\tilde q_1,q_2)|}{|n|}&\ll Q_1^5Q_2^{\frac{5}{2}}\sum_{1\leq |n|\leq X^{2013}}\frac{\sqrt{(n,q_2)}}{|n|}\ll Q_1^5Q_2^{\frac{5}{2}}X^{\varepsilon}.
\end{align}
Again applying Lemma \ref{char-sum-1}, it follows that
\begin{align}
\label{bd-char-2}
H\sum_{q_1\in\mathcal Q_1}\sum_{1\leq |n|\leq X^{2013}}&\frac{|\mathcal T(m,q_1n;q_1,q_1,q_2)|}{q_1|n|}+\frac{N}{QQ_1}\sum_{q_1\in\mathcal Q_1}|\mathcal T(m,0;q_1,q_1,q_2)|\\
\nonumber & \ll HQ^{5/2}X^{\varepsilon}+NQ^2X^{\varepsilon}.
\end{align}
The above two bounds \eqref{bd-char-1}, \eqref{bd-char-2} yield
\begin{align*}
\tilde D^\sharp(m,q_2)\ll (HQ_1^5Q_2^{5/2}+NQ^2)H(XY)^{2+\varepsilon}.
\end{align*}
Plugging this estimate in \eqref{dh} we get the following:
\begin{lemma}
\label{last-lemma}
For $Q_1Q_2=Q$, we have 
\begin{align*}
\tilde D_1(\alpha)\ll \frac{\sqrt{N}Q_2M}{Q^5}(\sqrt{H}Q_1^{5/4}Q^{5/4}+\sqrt{N} Q)\sqrt{H}YX^{1+\varepsilon}.
\end{align*}
\end{lemma}

The optimal breakup $Q_1Q_2=Q$ is now obtained by equating the two terms. We get that $Q_2=X^{2/5}$ and $Q_1=rX^{1/10+\delta}$. The optimal choice for $\delta$ is now obtained by equating the resulting error term with the previous error term, namely $X^{1-\delta}$. We get 
$$
\delta=\frac{1}{35}-\frac{2}{7}\frac{\log r}{\log X}.
$$ 
Finally one checks that the error term in \eqref{approx-pw-11} is satisfactory for the above choice of $\delta$. This holds as long as $r\ll X^{4/5}$.


\end{document}